\documentclass{amsart}
\usepackage{graphicx} 
\usepackage{amsmath}
\usepackage{amsthm}
\usepackage{amssymb}
\usepackage{comment}
\usepackage{nameref}
\usepackage{cite}
\usepackage{xcolor}
\numberwithin{equation}{section}

\newcommand{\ul}{U_-}
\newcommand{\ur}{U_+}

\newcommand{\de}{\delta}
\newcommand{\ep}{\epsilon}

\newcommand{\epw}{\delta_0}
\newcommand{\epe}{\epsilon_2}
\newcommand{\ept}{\epsilon_0}

\newcommand{\utl}{\Tilde{U}}

\newcommand{\x}{\mathbf{X}}

\newcommand{\sis}{\sigma}

\newcommand{\cfi}{c_f^{(i)}}

\newcommand{\si}{S_i^{\mathbf{X}_i}}

\newcommand{\dkk}{\;\partial_{x}{k_{S_i}^{\x_i}}}

\newcommand{\epss}{\delta_{S_i}}

\newcommand{\mn}{{\mu}}
\newcommand{\ro}{\mathbf{r_1}}

\newcommand{\ri}{\mathbf{r_i}}
\newcommand{\lo}{\mathbf{l_1}}

\newcommand{\li}{\mathbf{l_i}}

\newcommand{\stdy}{\mathbf{S}_i}
\newcommand{\sdy}{{\sigma}_{i}}
\newcommand{\uls}{U_-}
\newcommand{\urs}{U_+}

\newcommand{\wt}{\widetilde}
\newcommand{\s}{\sigma}

\newtheorem{theo}{Theorem}[section]
\newtheorem{lemm}[theo]{Lemma}

\newtheorem{assum}[theo]{Assumption}

\title[Asymptotic Stability of Weak Viscous Shocks]{$a$-contraction applied to the time-asymptotic stability of weak shocks for a general 1-D viscous system}
\date{\today}

\author{Young-Sam Kwon}
\address{Department of Mathematics, \newline Dong-A University, Busan 49315, Republic of Korea}
\email{ykwon@dau.ac.kr}

\author{Jiayun Meng}
\address{Department of Mathematics, \newline The University of Texas at Austin, Austin, TX 78712, USA}
\email{jiayun@utexas.edu}

\author{Alexis Vasseur}
\address{Department of Mathematics, \newline The University of Texas at Austin, Austin, TX 78712, USA}
\email{vasseur@math.utexas.edu}

\thanks{\textbf{Acknowledgment.}This work was partially funded by NSF-DMS 2306852, NSF-RTG 1840314, and NSF-EPSRC 2219434.}


\begin{document}

\begin{abstract}
    The method of $a$-contraction with shifts builds on the relative entropy method of Dafermos and DiPerna \cite{Darem,Direm} and allows both viscous shocks and rarefaction waves to be studied directly at the level of the solution. This unified framework led to the resolution of Matsumura's long-standing conjecture \cite{MNht} on the time-asymptotic stability of a viscous shock in composition with a rarefaction. Kang, Vasseur, and Wang \cite{MVWns} established this stability for the one-dimensional compressible barotropic Navier-Stokes system. Their result was subsequently extended to composite planar waves for general $n\times n$ viscous systems in arbitrary spatial dimension \cite{meng}.
    
    The present paper is a companion to that work. Here we prove the estimates on viscous shock profiles needed to apply $a$-contraction with shifts, including derivative estimates that go beyond the standard bounds. We also give a streamlined account of the main ideas of the method by focusing on a single viscous shock for a general one-dimensional system.
\end{abstract}
\maketitle
\tableofcontents
\section{Introduction}
The time-asymptotic stability of compositions of elementary waves, in particular of a viscous shock and a rarefaction, for viscous approximations of one-dimensional hyperbolic systems with genuinely nonlinear characteristic fields was first discussed by Matsumura and Nishihara in 1986 \cite{MNho}, and later posed as an open problem in 1992 \cite{MNht}. It remained unresolved until 2023 \cite{MVWns}. The main difficulty lies in the incompatibility between the classical methods used to study rarefactions and viscous shocks.

The energy method has long been a principal tool for studying the time-asymptotic stability of elementary waves. The spreading of rarefaction waves allows one to apply energy estimates directly at the level of the solution, whereas the compressive nature of viscous shocks makes such estimates difficult to close. Early proofs of shock stability \cite{MNvs,Gvs} therefore applied energy estimates to the antiderivative of the perturbation, an approach that requires the initial perturbation to have zero mass. To remove this restrictive assumption, the Green's function method \cite{Liu1985} and the Evans function method \cite{MZspectral,HLZspectral} have been developed. These approaches provide detailed pointwise estimates in one dimension, and the Evans function framework even allows the study of planar viscous shocks under multidimensional perturbations.

\smallskip

\noindent
\textbf{$a$-contraciton with shifts.} However, none of these classical approaches can be applied simultaneously to both rarefactions and shocks. As a result, the Matsumura conjecture remained open for decades until the introduction of the $a$-contraction with shifts method. The method of $a$-contraction with shifts is a recently developed energy method that builds on the relative entropy method introduced by Dafermos and DiPerna \cite{Darem,Direm}. By introducing a suitable weight \cite{KVa} and imposing an additional shift on the shock \cite{Leger}, this $L^2$-based method allows us to establish energy estimates directly at the level of the solutions and provides a unified framework for studying viscous shocks, rarefaction waves, and contact discontinuities. This framework has led to the stability of generic Riemann solutions for the compressible Navier-Stokes \cite{MVWns} and Navier-Stokes-Fourier equations \cite{KVWnsfgeneric}, resolving Matsumura's longstanding conjecture.

\noindent
\textbf{Multidimensional extensions.}  The first multidimensional extension of $a$-contraction with shifts established the stability of a single weak planar viscous shock for the three-dimensional Navier-Stokes system \cite{WWstabilitys}. Subsequent work extended this stability result to compositions of such shocks \cite{KLstabibilitys}. The result of Kang, Vasseur, and Wang \cite{MVWns} resolving Matsumura's conjecture has since been extended to general $n\times n$ viscous systems in arbitrary spatial dimension \cite{meng}. This extension establishes the stability of a composite wave consisting of a planar viscous shock and either a planar rarefaction or a second planar viscous shock.

\noindent
\textbf{Purpose of the paper.} The main mathematical contribution of this paper is to establish the estimates on viscous shock profiles needed to apply $a$-contraction with shifts. These include derivative estimates that go beyond the standard bounds. In the companion paper \cite{meng}, these estimates are used to complete the proof of time-asymptotic stability for composite planar waves in general multidimensional viscous systems. In the present paper, we also use the simpler setting of a single viscous shock in one space dimension to explain in detail how $a$-contraction with shifts yields the $L^2$ estimate underlying the stability proof.

\subsection{A general viscous system and its shock profiles}

We consider a general $n\times n$ $1$-D viscous system
\begin{equation} \label{eq for sol}
    \partial_t U + \partial_{x} f(U) = \partial_{x} \bigl(B(U)\partial_{x} \eta'(U)\bigl)\, ,  
\end{equation}
where $U: \mathbb{R}_{\geq 0}\times \mathbb{R}\rightarrow \mathbb{R}^n$. The flux function $f: \mathbb{R}^n\rightarrow \mathbb{R}^n$ is assumed to be smooth, strictly hyperbolic, and genuinely nonlinear. 
The viscosity coefficient matrix $B: \mathbb{R}^n \rightarrow \mathbb{R}^{n\times n}$ is assumed to be smooth and positive definite. The entropy $\eta:\mathbb{R}^n\to\mathbb{R}$ is assumed to be smooth and strictly convex, in the sense that its Hessian $\eta''$ is positive definite. Moreover, we assume that there exists an entropy flux $q: \mathbb{R}^n\rightarrow \mathbb{R}$ of $\eta$ such that 
\begin{equation}\label{def of entropy flux}
    \partial_i q(U) = \sum_{l=1}^n \partial_l \eta(U) \partial_i (f)_l (U) \, ,
\end{equation}
for any $1\leq i\leq n$. We endow the system with an initial value 
\begin{equation*}
    U(0, x) = U_0(x)\, , \quad x\in\mathbb{R} \, ,
\end{equation*}
with fixed end states $U_\pm \in \mathbb{R}^n$, i.e., 
\begin{equation} \label{initial value large time}
    U_0(x) \rightarrow U_\pm\, , \quad \textup{as } x\rightarrow \pm \infty \, .
\end{equation}
The long-time asymptotic behavior of solutions of \eqref{eq for sol} with initial value satisfying \eqref{initial value large time} is closely related to the one-dimensional Riemann problem for the corresponding inviscid model
\begin{equation}\label{eq for sol inviscid}
    \partial_t U + \partial_{x} f(U) = 0\, , 
\end{equation}
with the associated initial value 
\begin{equation}\label{inviscid initial value}
    U(0, x) = 
    \begin{cases}
         U_-  \, , \quad x<0 \, ,\\
         U_+  \, , \quad x>0 \, .
    \end{cases}
\end{equation}

\smallskip

\noindent
\textbf{Viscous shocks.}
In this paper, we focus on the case in which the Riemann initial data \eqref{inviscid initial value} give rise to a single shock in an extremal characteristic family. As descrbied in \cite{Dafermos2016}, the Rankine-Hugoniot set
\begin{equation}\label{Rankine-Hugoniot}
    \mathcal{S}(\ul) = \{U | -\sis (U-\ul) + f(U)-f(\ul) = 0 \text{ for a constant } \sis\}\, .
\end{equation}
consists, in a neighborhood of  $\ul$, of $n$ smooth curves $\mathcal{S}_1(\ul), ..., \mathcal{S}_n(\ul)$ intersecting at $\ul$. Wecall $\mathcal{S}_i(\ul)$ the $i$-shock curve. For any $\ur \in \mathcal{S}_i(\ul)$ sufficientlly close enough to $\ul$, there exists a weak solution $\Tilde{\mathbf{S}}_i$ to (\ref{eq for sol inviscid}) with initial value (\ref{inviscid initial value}) 
\begin{equation*}
    \Tilde{\mathbf{S}}_i(t,x) =
    \begin{cases}
        \ul \, , \quad x<\sis_i t \, , \\
        \ur \, , \quad x>\sis_i t \, , \\
    \end{cases}
\end{equation*}
where the Rankine-Hugoniot speed $\sis_i$.

For an admissible $i$-shock, viscosity in \eqref{eq for sol} smooths the jump discontinuity into a transition layer, called a viscous shock. We denote this traveling wave by $\mathbf{S}_i(x-\sis_i t)$, whose profile satisfies
\begin{equation}\label{eq for sol shock}
    -\sis_i \partial_{x} \mathbf{S}_i + \partial_{x} f(\mathbf{S}_i) = \partial_{x} \bigl(B(\mathbf{S}_i)\partial_{x} \eta'(\mathbf{S}_i)\bigl)\, .
\end{equation}

\subsection{Main results}

The following is the main result. We show that the long-time asymptotic state of a perturbation of a viscous $i$-shock wave is the viscous $i$-shock wave itself, up to a shift.

\begin{theo}\label{thm}
For any $\ul \in \mathbb{R}^n$, there exist constants $\epw, \ept>0$ such that the following is true.
\noindent
Let $i\in \{1,n\}$, and let $\ur \in \mathcal{S}_i(\ul)$ be such that
\begin{equation*}
    |\ur-\ul|< \epw \, .
\end{equation*}
We denote by $\mathbf{S}_i$ the viscous $i$-shock wave satisfying (\ref{eq for sol shock}) with end states $\ul$ and $\ur$. Let $U_0$ be an initial value such that 
\begin{equation}\label{thm initial data}
    \|U_0-U_-\|_{L^2(\mathbb{R}_{\leq 0})} + \|U_0-U_+\|_{L^2(\mathbb{R}_{\geq 0})} + \|\partial_x U_0\|_{L^2(\mathbb{R})}< \ept \, .
\end{equation}

\noindent
Then the viscous system (\ref{eq for sol}) has a unique global-in-time solution $U$. Moreover, there exist absolutely continuous shifts $\x_i(t)$ such that
\begin{align}
    &U(t,x) - \mathbf{S}_i\bigl(x-\sigma_i t + \x_i(t)\bigl) \in C\bigl(0,\infty;H^1(\mathbb{R})\bigl) \, , \label{thm continuity} \\
    &\underset{x\in \mathbb{R}}{\textup{sup}}\, \Bigl|U(t,x)-  \mathbf{S}_i\bigl(x-\sigma_i t + \x_i(t)\bigl)\Bigr| \rightarrow 0 \, , \quad \text{ as } \; t\rightarrow {{\infty}} \, . \label{thm pointwise convergence}
\end{align}
In addition,
\begin{align}
    U -\mathbf{S}_i\bigl(x-\sigma_i t + \x_i(t)\bigl) \in L^2\bigl(0,\infty;\dot{H}^{2}(\mathbb{R})\bigl) \, , \label{thm second deri}
\end{align}
and 
\begin{equation} \label{thm X shift}
    \underset{t\rightarrow \infty}{\textup{lim}} \dot{\mathbf{X}}_i(t)= 0 \, .
\end{equation}
\end{theo}
In this paper, $H^1$ denotes the Sobolev space $W^{1,2}$. We use $\dot{H}^m$ to denote the space of functions whose weak derivatives of orders $1,\ldots,m$ all belong to $L^2$.

\subsection{Structure of the paper}
One of the main purposes of this paper is to establish the properties of viscous shock profiles needed to apply $a$-contraction with shifts. We prove these properties in Appendix~\ref{sec proof S}. 

As in the companion paper, the proof of Theorem~\ref{thm} combines local-in-time estimates with a priori estimates. In the rest of this paper, we use the simpler setting of a single viscous shock in one space dimension to explain in detail how $a$-contraction with shifts yields the $L^2$ estimates underlying the a priori estimates. As shown in the companion paper \cite{meng}, the a priori estimates are established under the following assumption.
\begin{assum}\label{assumption}
    Let $U$ be solution to (\ref{eq for sol}) on $[0,T]$ for some $T>0$.
    Let $\si$ be the shifted viscous $i$-shock defined in (\ref{eq for super wave}) with absolutely continuous shifts $\x_i$ defined in (\ref{def of shift}) and weight function $a$ associated with it defined in (\ref{weight function}).
    Assume $\epss < \epw$,
    \begin{equation*}
        U - \si \in C\bigl([0,T]; H^1(\mathbb{R})\bigl) \, ,
    \end{equation*}
    and
    \begin{equation} \label{L infinity bound}
        \| U-\si\|_{L^\infty (0, T ; \, H^1 (\mathbb{R}))} \leq \epe \, .
    \end{equation}
\end{assum}
We discuss in detail the hyperbolic ``scalarization'' induced by the construction of the weight in Lemma \ref{lemma rem} and the Poincar\'e-type inequality enabled by the construction of the shift in Lemma \ref{lemma parabolic}. We collect the requirements from both steps and choose the strengths of the weight and the shift jointly to satisfy these requirements in Lemma \ref{lemma L^2 estimate}.

\section{Preliminaries}\label{sec preliminaries}

\subsection{Notations}
Before we go to the proofs, let us fix some notations. Let the eigenvalues of $f'(\ul)$ be $\lambda_1<...<\lambda_n$. For any $1\leq i\leq n$, let $\mathbf{r_i}$ and $\li$ be the right and left eigenvectors corresponding to $\lambda_i$ such that $\mathbf{r_i}$ is tangent to $\mathcal{S}_i(\ul)$, 
\begin{align}\label{def of li}
    \li=\eta''(\ul)\mathbf{r_i} \, , \quad \mathbf{r_i}\cdot \li=1 \, .
\end{align}
We define
\begin{equation} \label{r curve opposite s curve}
    \begin{aligned}
        \cfi :=( f''(\ul) : \mathbf{r_i} \otimes \mathbf{r_i}) \cdot \li &= \lambda_i '(\ul) \cdot \mathbf{r_i} = 2 \s_i'(\ul) \cdot \mathbf{r_i} < 0 \, .
    \end{aligned}
\end{equation}
Let $C$ denote the positive constants that depend only on $B, f, \eta, \ul$. Thanks to the existence of the entropy flux $q$, the matrix $\eta''f'$ is symmetric, and hence
\begin{equation} \label{trick eta'' f'}
    \li \cdot \mathbf{r_j} = \eta''(\ul) \mathbf{r_i} \cdot \mathbf{r_j} = 0 \text{ if } i\neq j \, .
\end{equation}

\subsection{Viscous shock wave}
Let $i\in \{1,n\}$. The viscous $i$-shock profile $\stdy$ satisfies
\begin{equation} \label{eq for viscous shock}
    \begin{cases}
        -\sdy \partial_{x} \stdy + \partial_{x} f(\stdy) = \partial_{x}\bigl(B(\stdy)\partial_{x} \eta'(\stdy)\bigl) \, , \\
        \stdy(-\infty) = \uls \, ,\; \stdy(+\infty) = \urs \, .
    \end{cases}
\end{equation}
Moreover, the Lax inequality implies
\begin{equation*}
    \lambda_i(\urs)< \sdy< \lambda_i(\uls) \, .
\end{equation*}
We define the wave strength of $\stdy$ by
\begin{align*}
    \epss = (\urs - \uls) \cdot \li \, .
\end{align*}
The existence of such profiles is established in \cite{MajdaSexist}. The following properties are direct consequences of Lemma \ref{lemma of shock}.
\begin{lemm}\label{pointwise bound of shock}
    For any $\uls \in \mathbb{R}^n$, there exist $\ep, C>0$ and $C_j>0$ for any $j\geq 2$ such that the following is true.

    \noindent
    For any $\urs \in \mathcal{S}_i(\uls)$ such that $|\uls-\urs| < \ep$, there exists a unique solution $\stdy$ to (\ref{eq for viscous shock}) such that
    \begin{align*}
        |\stdy({x})-\uls| &\leq C \epss e^{-C\epss|{x}|}, \quad \forall {x}< 0 \, , \\
        |\stdy({x})-\urs| &\leq C \epss e^{-C\epss|{x}|}, \quad \forall {x}> 0 \, , \\
        |\partial_{x} \stdy| &\leq C \epss^2 e^{-C\epss|{x}|}, \quad \forall {x} \in \mathbb{R}\, , \\
        |\partial_{x}^j \stdy| &\leq C_j \epss |\partial_{x} \stdy|, \quad \forall {x}\in \mathbb{R}\, , \; \forall j \geq 2 \, .
    \end{align*}
\end{lemm}
We call the following traveling profile the viscous $i$-shock wave
\begin{align*}
    S_i(t,x) = \stdy(x - \sigma_i t) \, .
\end{align*}
as the viscous $i$-shock wave. It satisfies 
\begin{align}\label{eq for shock}
    \partial_t S_i + \partial_{x} f(S_i) &= \partial_{x}\bigl(B(S_i)\partial_{x} \eta'(S_i)\bigl) \, .
\end{align}
We define the normalized variation of $S_i$ in the direction of $\ri$ by
\begin{equation}\label{decomposition of S}
    \begin{aligned}
        k_{S_i}(t,x) = \frac{\bigl(S_i(t,x) - \uls\bigl) \cdot \li}{\epss} \,.
    \end{aligned}
\end{equation}
By Lemma \ref{lemma of shock}, $k_{S_i}$ is strictly increasing, 
\begin{align}\label{derivative of k}
    \big|\partial_{x}k_{S_i} +\frac{\cfi}{2B(\ul)\li \cdot \li} \epss k_{S_i} (1-k_{S_i}) \big| \leq C \epss ^ {2} k_{S_i} (1-k_{S_i}) \, ,
\end{align}
and
\begin{equation}\label{k_S as main part}
    \begin{aligned}
        |\partial_{x} S_i-\epss \partial_{x} k_{S_i} \mathbf{r_i}| \leq C \epss^{2} \partial_{x} k_{S_i} \, .
    \end{aligned}
\end{equation}

\subsection{Construction of the weight functions and the shift function}\label{subsec construction}
We define the weight function $a_{S_i}$ by
\begin{equation*}
    \begin{aligned}
        a_{S_i} (t,x) = 
        \begin{cases}
            1 - \Lambda_{S_i} \epss{k_{S_i}}(t,x) \, , \quad \textup{if } i=1 \, , \\
            1 + \Lambda_{S_i} \epss  k_{S_i}(t,x)  \, , \quad \textup{if } i=n \, ,
        \end{cases}
    \end{aligned}
\end{equation*}
where the constant $\Lambda_{S_i}>0$ depends only on $\ul$ and is chosen sufficiently large in Lemma \ref{lemma L^2 estimate}. The weight function $a_{S_1}$ is decreasing and the weight function $a_{S_n}$ is increasing, since $k_{S_i}$ is increasing. Moreover, we take $\epss$ small enough so that 
\begin{align*}
    \frac{1}{2} \leq a_{S_i} \leq 2 \, .
\end{align*}enough.

We consider the shifted viscous $i$-shock wave 
\begin{equation} \label{eq for super wave}
    {S}_i^{{\mathbf{X}_i}}(t,x) = {S}_i\bigl(t,x+\x(t)\bigr) \, ,
\end{equation}
and the weight function associated with the shifted viscous $i$-shock 
\begin{align}\label{weight function}
    a(t,x) =  a_{S_i} (t,x+\x_i(t)) \, ,
\end{align}
where the shift $\mathbf{X}_i$ is defined as the solution to the ODE
\begin{equation} \label{def of shift}
    \begin{cases}
        \dot{\x}_i(t) = \frac{\Tilde{C}_i}{\epss} \int a \, \eta'' (\si)(U-\si) \partial_{x} S_i^{\x_i} \, dx\, , \\
        \x_i(0)=0 \,,
    \end{cases}
\end{equation}
where the constant $\Tilde{C_i}>0$ depends only on $\ul$ and is chosen sufficiently large in Lemma \ref{lemma L^2 estimate}. The existence and uniqueness of the shift $\x_i$ are established in \cite{meng}.

By (\ref{eq for shock}), we have
\begin{equation} \label{eq for superposition wave}
    \partial_t {S}_i^{{\mathbf{X}_i}} + \partial_{x} f({S}_i^{{\mathbf{X}_i}}) = \partial_{x} \bigl(B({S}_i^{{\mathbf{X}_i}})\partial_{x} \eta'  ({S}_i^{{\mathbf{X}_i}})\bigl)  + \dot{\mathbf{X}}_i\partial_{x}\si \, .
\end{equation}
Together with (\ref{eq for sol}) and (\ref{eq for superposition wave}), we get
\begin{equation}\label{eq u and u tilde}
    \begin{aligned}
        &\partial_t (U-\si) + \partial_{x} \bigl(f(U)-f(\si)\bigl) \\
        &= \partial_{x}\bigl(B(U)\partial_{x} \eta'(U) -B(\si)\partial_{x} \eta' (\si)\bigl) - \dot{\mathbf{X}}_i\partial_{x}\si \, .
    \end{aligned}
\end{equation}

\section{$L^2$ estimates}\label{section energy estimates}
We describe the main ideas behind the use of $a$-contraction with shifts to establish $L^2$ estimates for a general system. In Section \ref{subsec choice of basis}, we first construct a basis adapted to the flux structure for hyperbolic ``scalarization'' and to the viscosity structure for the Poincar\'e-type inequality. This basis contains the special direction associated with the viscous $i$-shock layer, and its vectors are mutually orthogonal with respect to the dissipation.

We then decompose the perturbation in this basis. In Section~\ref{subsec relative entropy}, we explain how the weight $a$ activates the spectral gap and generates favorable negative hyperbolic contributions in every direction except the special direction associated with the viscous shock layer. Their magnitude is proportional to $\Lambda_{S_i}$, the parameter defining the strength of the derivative of the weight.
Choosing $\Lambda_{S_i}$ sufficiently large allows these negative terms to compensate for all other hyperbolic contributions in those directions, leaving a positive scalar hyperbolic remainder in the special direction associated with the viscous shock layer.

The dissipation and the shift-induced term compensate for the positive hyperbolic remainder through a Poincar\'e-type inequality. More precisely, the orthogonality of the basis with respect to the dissipation ensures sufficient dissipation in the special
direction associated with the viscous layer to provide the necessary derivative control. We choose $\Tilde{C_i}$, the constant defining the strength of the shift, sufficiently large so that the shift-induced term provides the required control of the average. The details are given in Section \ref{subsec poincare}.

Finally, in Section \ref{subsec L^2 contraction}, we collect the requirements imposed by hyperbolic scalarization and the Poincar\'e-type inequality and choose $\Lambda_{S_i}$ and $\Tilde{C_i}$ jointly to satisfy them.

\subsection{Choice of basis}\label{subsec choice of basis}
As in \cite{meng}, we construct a basis that contains the special direction of the viscous shock layer and is orthogonal with respect to the dissipation. The first requirement comes from hyperbolic ``scalarization,'' in which we use the weight $a$ to exploit the spectral gap. The second is imposed to maximize the dissipation in the
special direction associated with the viscous shock layer, allowing us to control the hyperbolic remainder later through a Poincar\'e-type inequality.

Since the dissipation acts on entropy variables, we construct the basis in these variables, with $\li$ as the corresponding special direction associated with the viscous
shock layer. In particular, we choose a basis
\begin{align*}
    \{\li\} \cup \{\mathbf{v_j}\}_{j\neq i} \, ,
\end{align*}
and, for notational convenience, set $\mathbf{v}_i:=\li$, such that
\begin{align*}
    \bigl({B}(\ul) + B^T(\ul) \bigr) \mathbf{v_{j_1}} \cdot \mathbf{v_{j_2}} = 0 \, , \;
    \textup{for any } j_1 \neq j_2 \, , \quad {B}(\ul) \mathbf{v_{j}} \cdot \mathbf{v_{j}} = 1 \, , \; \textup{for any } j \neq i \, .
\end{align*}
This is always possible thanks to the Gram-Schmidt process. We decompose the perturbation in this basis:
\begin{equation} \label{decomposition of S_1}
    \eta'(U)(t,x) - \eta'(\si)(t,x) = \mu_i(t,x) \li + \sum_{j \neq i} \mu_j(t,x) \mathbf{v_j} \,.
\end{equation}
On the other hand, the hyperbolic analysis occurs at the level of the conserved quantity $U-\si$. By \eqref{L infinity bound}, we have
\begin{equation}\label{L infinity}
    \|U - \utl\|_{L^\infty((0,T)\times \mathbb{R})} \leq C \epe \, .
\end{equation}
Thanks to the Taylor expansion, we get
\begin{align}\label{decomposition of S_1 taylor}
    \big|(U-\si)-\bigl(\mu_i \ri + \sum_{j \neq i} \mu_j \, \eta''(\ul)^{-1}\mathbf{v_j}\bigl)\big| \leq C (\epe + \epw)|U-\si| \, .
\end{align}
Since $\eta''$ is positive definite,
\begin{align*}
    \{\ri\} \cup \{\eta''(\ul)^{-1} \mathbf{v_j}\}_{j\neq i}
\end{align*}
is also a basis. Hence, we have
\begin{align}\label{property of perturbation 1}
    c \sum_{j\neq i} |\mu_j|^2 \leq \bigl|\sum_{j\neq i} \mu_j \, \eta''(\ul)^{-1}\mathbf{v_j}\bigl|^2 \leq C \sum_{j \neq i} |\mu_j|^2 \, ,
\end{align}
for some constants $c,C>0$ that depend only on $\ul$.

We define $P$ to be the projection onto $\textup{span}\{\mathbf{r_j}\}_{j \neq i}$, i.e., 
\begin{align}\label{def of projection}
    P(v) = v - (v\cdot \li) \cdot \ri \, , \quad \textup{ for any } v\in \mathbb{R}^n \, .
\end{align}
We define the subspace
\begin{align}\label{def of V}
    V = \textup{span} \bigl\{\eta''(\ul)^{-1} \mathbf{v_j}\bigl\}_{j\neq i} \, .
\end{align}
We note that $\ri\notin V$, so there exists $C>0$ such that 
\begin{align}\label{property of perturbation 2}
    |P(v)|^2 \geq C |v|^2 \, \textup{ for any } v\in V\, .
\end{align}

\subsection{Hyperbolic ``scalarization"}\label{subsec relative entropy}
We study the evolution of the weighted relative entropy
\begin{equation*}
    \int a(t,x) \, \eta\bigl(U(t,x)|\si(t, x)\bigl) \,dx \, ,
\end{equation*}
where the relative entropy is defined by
\begin{equation}\label{eq for REF}
    \eta(U|V) = \eta(U) - \eta(V) - \eta'(V)(U-V) \, .
\end{equation}
We define the relative flux $f$
\begin{equation}\label{eq for relative flux}
    \begin{aligned}
        f(U|V) = f(U) - f(V) - f'(V)(U-V) \, ,
    \end{aligned} 
\end{equation}
and the relative entropy $q$ 
\begin{equation} \label{eq for flux of RE}
    q(U;V) = q(U) - q(V) - \eta'(V)\bigl(f(U) - f(V)\bigl) \, ,
\end{equation}
where $q(\cdot)$, the entropy flux of $\eta$, is defined in (\ref{def of entropy flux}).

In the following lemma, we decompose the perturbation in the basis introduced in Section \ref{subsec choice of basis} and show how this choice maximizes the dissipation in the special direction associated with the viscous layer. By analyzing the additional hyperbolic term $\mathcal{H}_1$ generated by the introduction of the weight, we also show how the weight activates the spectral gap and produces contributions with favorable negative signs in all directions except the special direction associated with the viscous shock layer.

The following lemma corresponds to Lemma~4.3 in \cite{meng}, but its proof is simpler in the present setting. Since the system is one-dimensional, there are no hyperbolic or dissipative terms in the transverse directions. Moreover, we consider only a single layer, so there are no errors arising from interactions between layers. We provide only a brief proof sketch focusing on the main ideas. For the remaining details, including the estimates for the shift-induced terms, we refer to the proof of the corresponding lemma in \cite{meng}.

As in \cite{meng}, we postpone our choice of $\Lambda_{S_i}, \Tilde{C}_i$ to subsection \ref{subsec L^2 contraction}, where we unify the choice of constants for both hyperbolic ``scalarization" and Poincar\'{e} inequality.
\begin{lemm} \label{lemma rem} 
    For any $\ul \in \mathbb{R}^n$, any $0<\beta,\delta_{\mathcal{H}}<1$, and any $\Lambda_{S_i}, \Tilde{C}_i>0$, there exist $\epw, \epe, C>0$ such that the following is true. 
    
    \noindent
    Assume Assumption \ref{assumption}. Then for any $t\in[0,T]$,
    \begin{equation*}
        \frac{d}{dt} \int a(t,x) \, \eta\bigl(U(t,x)|\si(t, x)\bigl) \,dx \leq \mathcal{Z}(U) - \mathcal{D}(U) + \mathcal{H}(U)    \, ,
    \end{equation*}
    where  
    \begin{equation*}
        \begin{aligned}
            \mathcal{Z}(U)&= \dot{\x}_i\mathcal{Y}_{i}(U) \,, \\
            \mathcal{D}(U) &= (1-\frac{\beta}{2}) \mathcal{D}^{p}(U) + C \mathcal{D}^{r}(U)  \, , \\
            \mathcal{H}(U) &= \bigl(-\mathbf{C}_1 \Lambda_{S_i} + \frac{\mathbf{C}_2}{\delta_{\mathcal{H}}}\bigl) \mathcal{H}_C + (1+\delta_{\mathcal{H}}) \mathcal{H}_{S_i}  \, , 
        \end{aligned}
    \end{equation*}
    where $\mathbf{C}_0, \mathbf{C}_1, \mathbf{C}_2 >0$ are constants that depend only on $\ul$, and
    \begin{align*}
        \mathcal{Y}_{i}(U) &= \int \partial_{x} a_{S_i} ^ {\mathbf{X}_i} \, \eta({U|\si})\, dx - \int a \eta''(\si)(U-\si) \partial_{x} S_i^{\mathbf{X}_i}\, dx \, ,\\
        \mathcal{D}^p(U) &= {B(\ul)\li \cdot \li}\int |\partial_{x}{\mu_i}| ^2 \, dx \, ,\\ 
        \mathcal{D}^{r}(U) &= \sum_{j\neq i} \int |\partial_{x} \mu_j | ^2 \, dx  \, ,\\ 
        \mathcal{H}_C(U) &= \sum_{j\neq i}\epss\int \mu_j ^ 2 \dkk dx \, ,  \\
        \mathcal{H}_{S_i}(U) &= -\frac{\cfi}{2} \epss \int \mu_i ^ 2 \dkk \, dx \, .
    \end{align*}
\end{lemm}

\noindent
\textit{Remark.}
We have
\begin{align*}
    \mathcal{D}^p \, , \; \mathcal{D}^r \, , \; \; \mathcal{H}_C \, , \; \mathcal{H}_{S_i} \geq 0 \, .
\end{align*}
We denote the shift-induced term by $\mathcal{Z}$ and the viscous dissipation by $\mathcal{D}$. Here, $\mathcal{D}^p$ corresponds to the special direction of the layer, while $\mathcal{D}^r$ corresponds to the remaining directions in the dissipation-orthogonal basis. Similarly, $\mathcal{H}$ denotes the hyperbolic contribution arising from the flux, with $\mathcal{H}_{S_i}$ corresponding to the special direction associated with the viscous shock layer and $\mathcal{H}_C$ to the remaining directions.

The construction of the weight yields $-\mathbf{C}_1\Lambda_{S_i}$ as one of the coefficients of $\mathcal{H}_C$ stated in the lemma. Choosing $\Lambda_{S_i}$ sufficiently large makes this coefficient negative, leaving a positive hyperbolic remainder only in the special direction of the layer. This is the hyperbolic ``scalarization.'' In Lemma~\ref{lemma parabolic} below, we control this remainder using $\mathcal{D}^p$ and the shift-induced term $\mathcal{Z}$ through a Poincar\'e-type inequality.

\begin{proof}
By the definition of relative entropy function (\ref{eq for REF}), (\ref{eq for sol}), and (\ref{eq for superposition wave}), we have
\begin{align*}
     &\frac{d}{dt} \int a(t,x) \, \eta\bigl(U(t,x)|\si(t, x)\bigl) \,dx \\
     &= \int \partial_t a \, \eta({U|\si})\, dx + \int a \,\Bigr[\bigl(\eta'(U) - \eta'(\si)\bigl)\partial_t U -  \eta'' (\si)(U-\si)\partial_t \si\Bigr ] dx \\ 
     &= \int \partial_t a \, \eta({U|\si})\, dx + \int a \, \biggr [\bigl(\eta'(U) - \eta'(\si)\bigl) \partial_{x}\bigl(B(U)\partial_{x} \eta'(U)\bigl) - \partial_{x} f(U) \\&\quad  
     - \eta''(\si)(U-\si)\Bigl(\partial_{x} \bigl(B(\si)\partial_{x} \eta' (\si)\bigl)  -\partial_{x} f(\si) + \dot{\mathbf{X}}_i\partial_{x}\si \Bigl)\biggr] \,dx \, .
\end{align*}
According to definitions of (\ref{eq for relative flux}) and (\ref{eq for flux of RE}), we have
\begin{equation*}
    \frac{d}{dt} \int a(t,x)\, \eta\bigl(U(t,x)|\si(t, x)\bigl) \, dx = \mathcal{Z}_0 + \sum_{i=1}^{3} \mathcal{D}_i + \sum_{i=1}^{2} \mathcal{H}_i\, ,
\end{equation*}
where the shift term is
\begin{equation*}
    \begin{aligned}
        \mathcal{Z}_0&= \int \partial_t a \, \eta({U|\si})\,dx  - 
            \dot{\mathbf{X}}_i\int a\, \eta''(\si)(U-\si) \partial_{x} \si \, dx \, , 
    \end{aligned}
\end{equation*}
the dissipation gives
\begin{align*}
    \mathcal{D}_1 &= \int a\, \bigl(\eta'(U)-\eta'(\si)\bigl)\partial_{x}\Bigl(B(U)\partial_{x} \bigl(\eta'(U)-\eta'(\si)\bigl)\Bigl) \,dx\, ,\\
    \mathcal{D}_2 &= \int a\, \bigl (\eta'(U)-\eta'(\si)\bigl)\partial_{x}\Bigl(\bigl(B(U)-B(\si)\bigl)\partial_{x} \eta'(\si)\Bigl) \,dx\, , \\
    \mathcal{D}_3 &= \int a \, \eta' (U|\si) \partial_{x}\bigl(B(\si)\partial_{x} \eta'(\si)\bigl) \,dx\, ,
\end{align*}
and the flux induces
\begin{align*}
    \mathcal{H}_1 &= \int \partial_{x} a \, q(U;\si)\,dx \, , \\
    \mathcal{H}_2 &= -\int a \, \partial_{x} \eta'(\si) f(U|\si) \,dx  \, . 
\end{align*}

\noindent
\textbf{Dissipation.} Following \cite{meng}, the main viscous term is
\begin{align*}
    \mathcal{D}_{11} = - \int a \, \partial_{x}\bigl(\eta'(U)-\eta'(\si)\bigl) B(U) \partial_{x} \bigl(\eta'(U)-\eta'(\si)\bigl) \,dx \, ,
\end{align*}
with the remaining viscous terms bounded by
\begin{align*}
    C(\epe+\epw^{1/4}) ( \mathcal{D}^{p} + 
    \mathcal{D}^{r})  + C \epw^{1/2}\int |\partial_{x} \si| |U-\si|^2 \, dx \, .
\end{align*}
We apply the dissipation-orthogonal decomposition of the perturbation given in \eqref{decomposition of S_1} to $\mathcal{D}_{11}$. By (\ref{L infinity}) and the fact that $B$ is positive definite, we take $\epe,\epw$ small enough and get
\begin{align*}
    \mathcal{D}_{11} &\leq -\bigl(1-C(\epe+\epw)-C \Lambda_{S_i}\epw\bigl) \\ 
    & \quad \int \bigl(\li \partial_{x}\mu_i + \sum_{j \neq i} \mathbf{v_j}\partial_{x}\mu_j\bigl)B(\ul) \bigl(\li \partial_{x}\mu_i + \sum_{j \neq i} \mathbf{v_j}\partial_{x}\mu_j\bigl)\, dx 
    \\&\leq - \bigl(1-C(\epe+\epw^{1/2})\bigl)\mathcal{D}^{p} - C \mathcal{D}^{r}\, .
\end{align*}

\noindent
\textbf{Hyperbolic ``scalarization".}
Following \cite{meng}, the main hyperbolic contributions from $\mathcal{H}_1$ and $\mathcal{H}_2$ are, respectively,
\begin{align*}
    I_1 &= J_1 + \mathcal{Z}_i \, , \\
    I_2 &= - \epss \int \dkk \,\li\cdot \frac{1}{2} f''(\ul) : (U-\si)\otimes (U-\si) \,dx \, ,
\end{align*}
where
\begin{align*}
    J_1 &= \frac{1}{2} \int \partial_{x} a (U-\si) \eta''(\ul) \bigl(f'(\ul)-\lambda_i \bigr)(U-\si) \, dx \, , \\
    \mathcal{Z}_i &= \frac{\lambda_i}{2} \int \partial_{x} a \, (U-\si) \eta''(\ul)(U-\si) \, dx \, .
\end{align*}
The remaining hyperbolic contribution is bounded by
\begin{align*}
    C(\epe^{1/2}+\epw^{1/2}) \int |\partial_{x} \si| |U-\si|^2 \, dx \, .
\end{align*}

The term $\mathcal{Z}_i$ arises from the propagation speed of the viscous shock layer and is analyzed together with the shift-induced term. We refer to \cite{meng} for the
details. Hyperbolic ``scalarization'' is made possible by the weight $a$ associated with the viscous shock, which generates the additional hyperbolic term $\mathcal{H}_1$. We apply to $J_1$ the decomposition in \eqref{decomposition of S_1 taylor}, which contains the special direction of the viscous shock layer. By (\ref{trick eta'' f'}), (\ref{decomposition of S_1 taylor}), and \eqref{def of projection}, we get
\begin{align}\label{scalarization}
    &J_1 - C (\epe^{1/2}+\epw^{1/2}) \int |\partial_{x} \si| |U-\si|^2 \, dx \\
    &\leq - \frac{\Lambda_{S_i}}{2} \min_{j\neq i} |\lambda_j-\lambda_i| \epss \int \eta'' (\ul)P(v) \cdot P(v) \dkk \, dx\, ,
\end{align}
where
\begin{align*}
    v = \sum_{j\neq i} \mu_j \eta''(\ul)^{-1} \mathbf{v_j} \in V  \, ,
\end{align*}
with $V$ defined in \eqref{def of V}. 

We note that obtaining the negative coefficient
\begin{align*}
    -\frac{\Lambda_{S_i}}{2} \min_{j\neq i} |\lambda_j-\lambda_i|
\end{align*}
requires the $i$-th characteristic family to be extremal, so that $\lambda_j-\lambda_i$ has the same sign for all $j \neq i$.

By (\ref{decomposition of S_1 taylor}), (\ref{property of perturbation 1}) and (\ref{property of perturbation 2}), we know
\begin{align*}
    P(v) \cdot P(v) \geq C \sum_{j\neq i} |\mu_j|^2 \, .
\end{align*}
Hence, we get
\begin{align*}
    &J_1 - C (\epe^{1/2}+\epw^{1/2}) \int |\partial_{x} \si| |U-\si|^2 \, dx \\
    &\leq - C {\Lambda_{S_i}} \min_{j\neq i} |\lambda_j-\lambda_i| \sum_{j \neq i} \epss  \int \mu_j^2 \dkk \, dx  \leq -C \Lambda_{S_i} \mathcal{H}_{C} \, .
\end{align*}

The main hyperbolic contribution from $\mathcal{H}_1$ provides cancellation in all basis directions except the special direction of the viscous shock layer. When estimating $I_2$ using the decomposition in \eqref{decomposition of S_1 taylor}, we therefore apply Young's inequality to make the hyperbolic contribution in this direction as small as possible:
\begin{equation*}
    \begin{aligned}
        I_2 -C(\epe+\epw)\int |\partial_{x} \si| |U-\si|^2\,dx 
        \leq (1+\frac{\delta_{\mathcal{H}}}{2}) \mathcal{H}_{S_i} + \frac{C}{\delta_{\mathcal{H}}} \mathcal{H}_C  \, .
    \end{aligned}
\end{equation*}
In all, we get
\begin{equation*}
    \begin{aligned}
        &\mathcal{H}_1 + \mathcal{H}_2 - C (\epe^{1/2}+\epw^{1/2}) \int |\partial_{x} \si| |U-\si|^2 \, dx - \mathcal{Z}_i \\ 
        &\leq \bigl(-C\Lambda_{S_i} + \frac{C}{\delta_{\mathcal{H}}}\bigl) \mathcal{H}_C + (1+\frac{\delta_{\mathcal{H}}}{2}) \mathcal{H}_{S_i} \, .
    \end{aligned}
\end{equation*}

\end{proof}
\subsection{Poincar\'{e} type inequality}\label{subsec poincare} 
The following Poincar\'e-type inequality was proved in \cite{KVpoincare}.
\begin{lemm}\label{Poincare inequality lemma}
    For any $f:[0,1]\rightarrow\mathbb{R}$ satisfying $\int_0^1 y (1-y)|f'|^2 \, dy < \infty$,
    \begin{equation}
        \int_0^1 |f-\int_0^1 f \, dy|^2 \, dy \leq \frac{1}{2} \int_0^1 y (1-y)|f'|^2 \, dy \, .
    \end{equation}
\end{lemm}
In the following lemma, we reduce the compensation of the positive hyperbolic remainder $\mathcal{H}_{S_i}$ by the dissipation $\mathcal{D}^p$ in the special direction and the shift-induced term $\dot{\x}_i\mathcal{Y}_i$ to an application of this Poincar\'e-type inequality.

By choosing $\Tilde{C}_i$, the constant defining the strength of the shift in \eqref{def of shift}, sufficiently large, we ensure that the only positive remainder in the lemma's inequality is the hyperbolic contribution $\mathcal{H}_C$. Hyperbolic ``scalarization'' allows us to compensate for this remainder by choosing $\Lambda_{S_i}$ sufficiently large, as shown in Lemma \ref{lemma L^2 estimate} in the next section.

\begin{lemm}\label{lemma parabolic}
    For any $\ul \in \mathbb{R}^n$, any $0<\beta<\frac{1}{3}$, and any $\Lambda_{S_i}, \Tilde{C}_i>0$, there exist $\epw, \epe, C>0$ such that the following is true. 
    
    \noindent 
    Assume Assumption \ref{assumption}. Then for any $t\in[0,T]$,
    \begin{align*}
    &\mathcal{H}_{S_i}(U) - (1-\beta) D^{p}(U) + \dot{\x}_i(t) \mathcal{Y}_i(U) \\ 
    &\leq -C \int |\partial_{x}{\mn_i}| ^2 \, dx- C \epss \int \mn_i ^ 2 \dkk \, dx- \frac{\epss}{2\Tilde{C}_i}|\dot{\mathbf{X}}_i(t)| ^2  \\
    &\quad + \bigl(- \frac{\mathbf{C}_3}{\epss} \Tilde{C}_i  + \frac{\mathbf{C}_4}{\epss} \bigl)\bigl(\epss \int \mn_i \dkk \, dx \bigl)^2
    + \mathbf{C}_5 \Tilde{C}_i \mathcal{H}_C(U) \, ,
    \end{align*}
    where $\mathbf{C}_3, \mathbf{C}_4, \mathbf{C}_5>0$ are constants that depend only on $\ul$, $\x_i$ is defined in \eqref{def of shift}, and $\mathcal{H}_{S_i},\mathcal{H}_C, \mathcal{D}^p, \mathcal{Y}_i$ are defined in Lemma \ref{lemma rem}.
\end{lemm}
\begin{proof}
First, we fix the time $t$. As ${k_{S_i}^{\x_i}}$ is strictly increasing, we define the change of variable
\begin{equation*}
    x\mapsto z:= {k_{S_i}^{\x_i}}(t,x) \, , \quad h(z) = \mn_i\bigl(t,x\bigl) \, .
\end{equation*}
By (\ref{derivative of k}), we can take $\epw$ small enough and get
\begin{align*}
    \mathcal{H}_{S_i} &= - \frac{\cfi}{2} \epss \int_0^1 h(z) ^ 2  \, dz  \, , \\
    -(1-\beta)D^{p} &\leq (\frac{1}{2}-\frac{3\beta}{4}) \, \cfi \epss \int_0^1 h'(z)^2 z (1-z) \,dz \, .
\end{align*}
Lemma \ref{Poincare inequality lemma} gives
\begin{align*}
    \int_0^1 h(z)^2\, dz &= \int_0^1 \bigl(h(z) - \int_0^1 h(z)\,dz\bigl)^2  \, dz +  \bigl(\int_0^1 h(z)\,dz\bigl) ^ 2\\
    &\leq \frac{1}{2} \int_0^1 h'(z)^2 z(1-z)\, dz + \bigl(\int_0^1 h(z)\,dz\bigl) ^ 2\, .
\end{align*}
Since $0<\beta<\frac{1}{3}$, we get
\begin{equation}\label{poincare inequality}
    \begin{aligned}
        &\mathcal{H}_{S_i}-(1-\beta)D^{p} \\ 
        &\leq - C \epss \int \mn_i ^ 2 \dkk \, d x -C \int |\partial_{x}{\mn_i}| ^2 \, d x  + \frac{C}{\epss} \bigl(\epss \int \mn_i \dkk \, d x\bigl)^2 \, .
    \end{aligned}
\end{equation}
It remains to verify that the shift-induced term $\dot{\x}_i\mathcal{Y}_i$ obtained in Lemma~\ref{lemma rem} controls the last positive term, namely the square of the
average. Recall the definition of the shift $\x_i$ in \eqref{def of shift}. Following the proof of Lemma~4.4 in \cite{meng}, we have
\begin{align*}
    &\Bigl|\dot{\x}_i\mathcal{Y}_i  + \frac{\epss}{\Tilde{C}_i}|\dot{\mathbf{X}}_i(t)|^2 \Bigr| + 
    \Bigl| \frac{\epss}{\Tilde{C}_i}|\dot{\mathbf{X}}_i(t)|^2
    - \frac{\Tilde{C}_i}{\epss}\bigl(\int  \bigl(\eta'(U) - \eta'(\si)\bigl)  \partial_{x} \si  \, dx\bigl)^2 \Bigr| \\
    &\leq \frac{C}{\epss}(\epe^{1/2}+\epw^{1/2}) \bigl(\int |U-\si| |\partial_{x_1} \si| \, dx\bigl)^2 \, .
\end{align*}
By  (\ref{k_S as main part}), Young's inequality, and H\"{o}lder's inequality, we take $\epe,\epw$ small enough and get
\begin{align*}
    &\bigl(\int  \bigl(\eta'(U) - \eta'(\si)\bigl)  \partial_{x} \si  \, dx\bigl)^2 \\
    &\geq C \bigl(\epss \int \mn_i \dkk \, dx\bigl)^2 -  C\sum_{j\neq i} \bigl(\epss \int \mn_j \dkk \, dx\bigl)^2 \\ &\quad-C (\epe^{1/2}+\epw^{1/2}) \,\epss^2 \int |U-\si|^2 \dkk\, dx \\
    &\geq  C \bigl(\epss \int \mn_i \dkk \, dx \bigl)^2 - C \sum_{j\neq i}\epss ^2 \int \mn_j ^2 \dkk \, dx \\
    &\quad- C (\epe^{1/2}+\epw^{1/2}) \, \epss^2 \int \mn_i ^2 \dkk \, dx \, .
\end{align*}
Hence, we take $\epe, \epw$ small enough and get
\begin{equation*}
    \begin{aligned}
        &\dot{\x}_i\mathcal{Y}_i + \frac{\epss}{2\Tilde{C}_i}|\dot{\mathbf{X}}_i(t)| ^2 \\
        &\leq - \frac{C}{\epss} \Tilde{C}_i \bigl(\epss \int \mn_i \dkk \, dx \bigl)^2 + C \Tilde{C}_i \sum_{j\neq i} \epss \int \mn_j ^2 \dkk \, dx  \\&\quad + C (\epe^{1/4}+\epw^{1/4}) \, \epss \int \mn_i ^2 \dkk \, dx \, .
    \end{aligned}
\end{equation*}
Together with \eqref{poincare inequality}, we take $\epe,\epw$ small enough and get the estimate.
\end{proof}
\subsection{$L^2$ estimates}\label{subsec L^2 contraction}
We now choose the constants $\Lambda_{S_i}$ and $\Tilde{C}_i$, which define the weight and the shift, respectively, in a way that work for both Lemma~\ref{lemma rem} and Lemma~\ref{lemma parabolic}.
\begin{lemm} \label{lemma L^2 estimate}
    For any $\ul \in \mathbb{R}^n$, there exist $\epw, \epe, C, \Lambda_{S_i}, \Tilde{C}_i > 0$ such that the following is true. 

    \noindent
    Assume Assumption \ref{assumption}. Then for any $t \in [0,T]$,
     \begin{align*}
        \begin{split}
            \|U(t,\cdot)-\si(t,\cdot)\|_{L^2(\mathbb{R})} ^ 2
            + \int_0^t D_0(U) + G_{S_i}(U)  + Y \, ds 
            \leq  C \|U_0-\si(0,\cdot)\|_{L^2(\mathbb{R})} ^2 \, ,
        \end{split}
     \end{align*}
    where 
    \begin{align*}
        D_0 &= \int \big|\partial_{x} (U-\si)\big|^2 \, dx \, , \\
        G_{S_i}(U)&=\int |U-\si|^2 |\partial_{x} \si| \, dx\, , \\
        Y &= \epss |\dot{\x}_i|^2 \, .
    \end{align*}
\end{lemm}
\begin{proof}
We first fix $0<\beta< \frac{1}{3}$ and then choose $0<\delta_{\mathcal{H}}<1$ small enough such that
\begin{align*}
    (1+\delta_{\mathcal{H}})(1-\beta) < 1-\frac{\beta}{2} \, .
\end{align*}
We take
\begin{align*}
    \Tilde{C}_i = \frac{2\mathbf{C}_4}{\mathbf{C}_3} \, .
\end{align*}
We take  $\Lambda_{S_i}$ large enough such that
\begin{align*}
    -\mathbf{C}_1 \Lambda_{S_i} + \frac{\mathbf{C}_2}{\delta_{\mathcal{H}}} + (1+\delta_{\mathcal{H}})\mathbf{C}_5\Tilde{C}_i =  -1 \, .
\end{align*}
Finally, we take $\epe, \epw$ small enough such that Lemma \ref{lemma rem} and Lemma \ref{lemma parabolic} give the estimate.
\end{proof}

\appendix
\section{Proof of Lemma \ref{pointwise bound of shock}}\label{sec proof S}
In this section, we prove Lemma~\ref{pointwise bound of shock} for $i=1$. The case $i=n$ is analogous. To simplify notation, we denote the shock satisfying \eqref{eq for viscous shock} by $S$, its speed by $\s$, and its left and right states by $U_-$ and $U_+$, respectively. 

We first define 
\begin{align*}
    \de = (U_+-U_-) \cdot \mathbf{l_1} \, , \quad k(x) = \frac{\bigl(S(x)-U_-\bigr) \cdot \mathbf{l_1}}{\de} \, .
\end{align*}
Here, $k$ measures the normalized variation of $S$ in the direction of $\ro$. By \cite{MajdaSexist}, we have $0<k(x)<1$, and there exists $C>0$ such that
\begin{align*}
    |S-U_-| \leq C \de \, .
\end{align*}
Decomposing with respect to the basis of right eigenvectors of $f'(U_-)$, we get
\begin{equation}\label{proof of S decomposition}
    \begin{aligned}
        U_+ - U_- = \de \mathbf{r_1} + \de^2 \sum_{i=2}^n c_i \mathbf{r_i} \, , \\
        S - U_- = \de k \mathbf{r_1} + \de^2 \sum_{i=2}^n g_i \mathbf{r_i} \, .
    \end{aligned}
\end{equation}
for some constants $c_i$ and some functions $g_i$.

The following lemma shows that the variation of $S$ in the direction of $\ro$ dominates its variation in the remaining directions. More precisely, $k$ controls $g_i$, and $k'$ controls $g_i'$. We also obtain an explicit formula for the first derivative that describes its decay as $|x|\to\infty$. Finally, we show that the first derivative controls all higher-order derivatives. These estimates directly imply Lemma \ref{pointwise bound of shock}.
\begin{lemm}\label{lemma of shock}
    For any $U_- \in \mathbb{R}^n$, there exist $\epsilon, C, C_j>0$ such that the following is true.
    Assume $|U_+-U_-|<\epsilon$. Then
    \begin{itemize}
        \item[1)] $\sum_{i=2}^n |g_i| \leq C |k|$ and $\sum_{i=2}^n |c_i - g_i| \leq C |1-k|$.
        \item[2)] $\bigl|B(S) \eta''(S) \partial_x S \cdot \mathbf{l_1} + \frac{ c_f^{(1)}}{2} \de^2 k (1-k)\bigr| \leq C \de^3 k (1-k)$.
        \item[3)] $\sum_{i=2}^n |g_i'| \leq C|k'|$.
        \item[4)] For any $j\geq 2$, $|\partial_{x}^j S|\leq C_j \de |\partial_x S|$.
    \end{itemize}
\end{lemm}

\subsection{Preliminary}
We begin by proving the following lemma, in which we use Taylor expansions to express $k'$ and the $g_i'$ in terms of $k$ and the $g_i$.
\begin{lemm}
    We have
    \begin{align}\label{proof of S 1) k and g}
        \begin{aligned}
            &|g_i' - A_i^L| \leq C \de |k| + C \de \sum_{j=2}^n |g_j| \, , \quad
            |g_i' - A_i^R| \leq C \de |1-k| + C \de \sum_{j=2}^n |c_j - g_j| \, , \\
            &|k'| \leq C \de|k| + C\de \sum_{i=2}^n|g_i|\, , \quad
            |k'| \leq C \de|1-k| + C\de \sum_{i=2}^n|c_i - g_i| \, ,
        \end{aligned}
    \end{align}
    where
    \begin{equation}\label{proof of S 3) A_L, A_R}
        \begin{aligned}
            A_i^L &= (-k+k^2) \frac{c_f^{(1)}}{2} d_{i1} + k^2 C_i + \sum_{j=2}^n g_j\bigl(\lambda_j-\lambda_1\bigl) d_{ij} \, , \\
            A_i^R &= (-k+k^2) \frac{c_f^{(1)}}{2} d_{i1} - (1-k^2) C_i - \sum_{j=2}^n (c_j-g_j)\bigl(\lambda_j-\lambda_1\bigl) d_{ij} \, ,
        \end{aligned}
    \end{equation}    
    with
    \begin{align*}
        C_i = \sum_{j=2}^n \frac{1}{2} f''(U_-):\mathbf{r_1}\otimes\mathbf{r_1} \cdot \mathbf{l_j} \, d_{ij} \, .
    \end{align*}
\end{lemm}
\begin{proof}
Integrating \eqref{eq for viscous shock} from $\pm\infty$ to $x$, we get
\begin{equation}\label{eq for viscous shock integrated}
    B_S := B(S)\eta''(S) \partial_x S 
    \begin{cases}
        =-\s(S-U_-) + f(S) - f(U_-) =: F_L \, , \\
        =-\s(S-U_+) + f(S) - f(U_+) =: F_R \, .
    \end{cases}
\end{equation}
Taylor expansion gives
\begin{equation*}\label{proof of S Taylor expansion}
    \begin{aligned}
        \bigl| B_S - B(\ul)\eta''(\ul) \partial_x S \bigr| &\leq C \de |\partial_x S| \, , \\
        \Bigl| F_L - \bigl(f'(\ul)-\s\bigr)(S-\ul) - \frac{1}{2} f''(\ul):(S-\ul)\otimes(S-\ul) \Bigr| &\leq C |S-\ul|^3 \, , \\
        \Bigl| F_R - \bigl(f'(\ur)-\s\bigr)(S-\ur) - \frac{1}{2} f''(\ur) : (S - \ur) \otimes (S-\ur) \Bigr| &\leq C |S-\ur|^3 \, .
    \end{aligned}
\end{equation*}
The last inequality further gives
\begin{align*}
    &\Bigl| F_R - \bigl(f'(\ul)-\s\bigr)(S-\ur) - f''(\ul) : (\ur - \ul) \otimes (S-\ur) \\
    &\qquad - \frac{1}{2} f''(\ul): (S-\ur) \otimes (S-\ur) \Bigr| \leq C \de^2 |S-\ur| \, .
\end{align*}
We next substitute the decomposition in \eqref{proof of S decomposition} into the above Taylor expansions of $B_S$, $F_L$, and $F_R$. We note first that \eqref{def of li} implies
\begin{align*}
    B(\ul)\eta''(\ul) \partial_x S  = \de k' B(U_-)\mathbf{l_1} + \de^2 \sum_{i=2}^n g_i' B(U_-)\mathbf{l_i} \, ,
\end{align*}
and \eqref{r curve opposite s curve} gives
\begin{align}\label{proof of S prem lambda sigma}
    \Bigl|\bigl(\lambda_1-\s\bigr) \mathbf{r_1} + \de \frac{c_f^{(1)}}{2} \mathbf{r_1} \Bigr| \leq C \de^2 \, . 
\end{align}
Hence, we get
\begin{equation}\label{eq for viscous shock integrated each side}
    \begin{aligned}
        &|B_S - \de k' B(U_-)\mathbf{l_1} - \de^2 \sum_{i=2}^n g_i' B(U_-)\mathbf{l_i}| \leq C \de^2 |k'| + C \de^3 \sum_{i=2}^n |g_i'| \, , \\
        &|F_L + \de^2 k \frac{c_f^{(1)}}{2} \mathbf{r_1} - \de^2 \sum_{i=2}^n g_i \bigl(\lambda_i-\lambda_1\bigl) \mathbf{r_i} \\
        &\qquad - \de^2 k^2 \frac{1}{2}f''(U_-): \mathbf{r_1}\otimes \mathbf{r_1}| \leq C \de^3 |k| + C \de^3 \sum_{i=2}^n |g_i| \, ,  \\
        &|F_R - \de^2 (1-k) \frac{c_f^{(1)}}{2} \mathbf{r_1} + \de^2 \sum_{i=2}^n (c_i-g_i) \bigl(\lambda_i-\lambda_1\bigl) \mathbf{r_i} \\
        &\qquad + \de^2 (1-k^2) \frac{1}{2}f''(U_-): \mathbf{r_1}\otimes \mathbf{r_1}| \leq C \de^3 |1-k| + C \de^3 \sum_{i=2}^n |c_i - g_i| \, .
    \end{aligned}
\end{equation}
We first take the scalar product of both sides of \eqref{eq for viscous shock integrated} with $\mathbf{l_1}$. By \eqref{eq for viscous shock integrated each side}, we get
\begin{align*}
    &|B_S \cdot \lo - \de k' B(U_-)\lo \cdot \lo| \leq C \de^2 |k'| + C \de^2 \sum_{i=2}^n |g_i'| \, , \\
    &|F_L \cdot \lo| \leq C \de^2 |k| + C \de^3 \sum_{i=2}^n |g_i| \, , \quad
    |F_R \cdot \lo| \leq C \de^2 |1-k| +  C \de^3 \sum_{i=2}^n |c_i - g_i| \, ,
\end{align*}
since \eqref{trick eta'' f'} implies
\begin{align}\label{proof of S prem scalar l1}
    \de^2 \sum_{i=2}^n g_i \bigl(\lambda_i-\lambda_1\bigl) \mathbf{r_i} \cdot \lo = 0 = \de^2 \sum_{i=2}^n (c_i-g_i) \bigl(\lambda_i-\lambda_1\bigl) \mathbf{r_i} \, .
\end{align}
Together with the fact that $B$ is positive definite, we get
\begin{align}\label{proof of S 1) k}
    \begin{aligned}
        |k'| &\leq C \de\sum_{i=2}^n|g_i'| + C \de|k| + C\de^2 \sum_{i=2}^n|g_i|\, , \\
        |k'| &\leq C \de\sum_{i=2}^n|g_i'| + C \de|1-k| + C\de^2 \sum_{i=2}^n|c_i - g_i|\, .
    \end{aligned}
\end{align}
On the other hand, for any $2\leq i\leq n$, there exist constants
$d_{ij}$ such that
\begin{align}\label{proof of S 3) d def}
    (B^{-1})^T(U_-) \mathbf{r_i} = \sum_{j=1}^n d_{ij} \mathbf{l_j} \, .
\end{align}
We also take the scalar product of both sides of \eqref{eq for viscous shock integrated} with $(B^{-1})^T(U_-) \mathbf{r_i}$. By \eqref{trick eta'' f'}, we note that
\begin{align*}
    \de k' B(U_-)\mathbf{l_1} \cdot (B^{-1})^T(U_-) \ri = 0 \, , \quad \Bigl( \de^2 \sum_{j=2}^n g_j' B(U_-)\mathbf{l_j} \Bigr) \cdot (B^{-1})^T(U_-) \ri = \de^2 g_i' \, .
\end{align*}
Hence, \eqref{eq for viscous shock integrated each side} give
\begin{align*}
    |B_S \cdot (B^{-1})^T(U_-) \ri - \de^2 g_i'| 
    \leq C\de^2 |k'| + C \de^3 \sum_{j=2}^n |g_j'| \, .
\end{align*}
Together with \eqref{proof of S 1) k}, we get
\begin{equation}\label{proof of S prelim scalar product li B}
    \begin{aligned}
        |B_S \cdot (B^{-1})^T(U_-) \ri - \de^2 g_i'| &\leq C \de^3 \sum_{j=2}^n |g_j'| + C \de^3 |k| + C\de^4 \sum_{i=2}^n|g_i| \, , \\
        |B_S \cdot (B^{-1})^T(U_-) \ri - \de^2 g_i'| &\leq C \de^3 \sum_{j=2}^n |g_j'| + C \de^3 |1-k| + C\de^4 \sum_{i=2}^n|c_i - g_i| \, .
    \end{aligned}
\end{equation}
On the other hand, \eqref{eq for viscous shock integrated each side} and \eqref{proof of S prem lambda sigma} give
\begin{equation}\label{proof of S prelim scalar product li F}
    \begin{aligned}
        |F_L - \de^2 A_i^L| &\leq C \de^3 |k| + C \de^3 \sum_{i=2}^n |g_i| \, , \\
        |F_R - \de^2 A_i^R| &\leq C \de^3 |1-k| + C \de^3 \sum_{i=2}^n |c_i - g_i| \, , 
    \end{aligned}
\end{equation}
where $A_i^L, A_i^R$ are defined in \eqref{proof of S 3) A_L, A_R}. Combining \eqref{proof of S prelim scalar product li B} and \eqref{proof of S prelim scalar product li F}, we conclude 
\begin{align}\label{proof of S 1) g pre}
    \begin{aligned}
        |g_i' - A_i^L| &\leq C \de \sum_{j=2}^n |g_j'| + C \de |k| + C \de \sum_{j=2}^n |g_j| \, , \\
        |g_i' - A_i^R| &\leq C \de \sum_{j=2}^n |g_j'| + C \de |1-k| + C \de \sum_{j=2}^n |c_j - g_j| \, .
    \end{aligned}
\end{align}
Since 
\begin{align*}
    |A_i^L| \leq C|k| + C\sum_{j=2}^n |g_j| \, , \quad |A_i^R| \leq C |1-k| +  C\sum_{j=2}^n |c_j-g_j| \, ,
\end{align*}
we conclude \eqref{proof of S 1) k and g} from \eqref{proof of S 1) g pre} and \eqref{proof of S 1) k}.

\end{proof}

\subsection{Proof of 1)}
\subsubsection{Proof of $|g_i|\leq C|{k}|$}\label{subsub proof of S 1) L}
By \eqref{proof of S 1) k and g} and \eqref{proof of S 3) A_L, A_R}, we get
\begin{align}\label{proof of S 1) k,g L}
    \begin{aligned}
        |k'| \leq C \de|k| + C\de\sum_{i=2}^n|g_i|\, ,\quad
        |g_i' - \sum_{j=2}^n g_j\bigl(\lambda_j-\lambda_1\bigl) d_{ij}| \leq C |k| + C \de \sum_{j=2}^n |g_j| \, .
    \end{aligned}
\end{align}
We define
\begin{equation}\label{def of G, D, Lamda}
    \begin{aligned}
        G = \begin{pmatrix}
            g_2 \\ \vdots \\ g_n 
        \end{pmatrix} , \quad 
        D = \begin{pmatrix}
            d_{22} & \cdots & d_{2n}\\
            \vdots & \ddots & \vdots\\
            d_{n2} & \cdots & d_{nn}
        \end{pmatrix}, \\
        \Lambda = \begin{pmatrix}
            \lambda_2-\lambda_1 & 0 & \cdots & 0\\
            0 & \lambda_3 -\lambda_1 & \cdots & 0\\
            \vdots & \vdots & \ddots & \vdots\\
            0 & 0 & \cdots & \lambda_n -\lambda_1
        \end{pmatrix}.
    \end{aligned}
\end{equation}
We note that \eqref{proof of S 3) d def} and \eqref{trick eta'' f'} imply
\begin{align*}
    D = R^T (B^{-1})(U_-) R \, ,
\end{align*}
where $R = \begin{pmatrix} \mathbf{r_2} \cdots \mathbf{r_n} \end{pmatrix}$. We get that $D$ is positive definite. 

We let
\begin{align*}
    h = \frac{G \cdot \Lambda G}{2{k}^2} \, .
\end{align*}
Since $\Lambda$ is positive definite, it suffices to show that $|h|$ is uniformly bounded on $\mathbb{R}$ to complete the proof. To establish this bound, we compute the derivative of $h$:
\begin{align*}
    h' = \frac{G'\cdot \Lambda G}{{k}^2} - 2\frac{{k}'}{{k}} h \, .
\end{align*}
By \eqref{proof of S 1) k,g L}, we get
\begin{align*}
    |G' - D \Lambda G| \leq C|k| + C \de |G| \, ,
\end{align*}
and thus
\begin{align}\label{proof of S 1) L h}
    |h' - \mathbf{C}_D h| \leq C_{1/2} h^{1/2} + C_1 {\de} h + C_{3/2} {\de} |h|^{3/2} \, ,
\end{align}
where $C_{1/2},C_1,C_{3/2}>0$ are constants depending only on $U_-$, and $\mathbf{C}_D(x)$ is a function. Since $D$ is positive definite, there exist $c_D, C_D>0$, depending only on $D,\Lambda$, such that
\begin{align*}
    c_D\leq \mathbf{C}_D \leq C_D \, .
\end{align*}
It suffices to show $h$ is bounded above by a constant on $\mathbb{R}$. We note as 
\begin{align*}
    \lim_{x\rightarrow \infty} g_i(x) = c_i \, , \quad \lim_{x\rightarrow \infty} {k}(x) = 1 \, ,
\end{align*}
there exist $x_* \in \mathbb{R}$ and $\wt{C}>0$ such that
\begin{align}\label{proof of S 1) contradiction L}
    h \leq \wt{C} \, , \quad \textup{on } [x_*,\infty) \, .
\end{align}
We take $M>\wt{C}$ to be big enough such that
\begin{align}\label{proof of S def M}
    \frac{c_D M}{2} > C_{1/2} M^{1/2} \, .
\end{align}
Then we can take $\de>0$ small enough such that 
\begin{align}\label{proof of S property M}
    c_D M - C_{1/2} M^{1/2} - C_1 {\de} M - C_{3/2} {\de} M^{3/2} > \frac{c_D M}{4} \, .
\end{align}
By \eqref{proof of S 1) L h}, whenever $h=M$, then
\begin{align*}
    h' > \frac{c_D M}{4} \, .
\end{align*}
Thus, if $h(x_{**})=M$ at some point $x_{**}$, then $h\geq M$ on $[x_{**},\infty)$. This contradicts \eqref{proof of S 1) contradiction L}. Therefore, $h\leq M$ on $\mathbb{R}$.

\subsubsection{Proof of $|c_i - g_i|\leq |1-{k}|$}
By \eqref{proof of S 1) k and g} and \eqref{proof of S 3) A_L, A_R}, we get
\begin{align}\label{proof of S 1) k,g R}
    \begin{aligned}
        |k'| &\leq C \de|1-k| + C\de\sum_{i=2}^n|c_i - g_i|\, , \\
        |g_i' + \sum_{j=2}^n (c_j-g_j)\bigl(\lambda_j-\lambda_1\bigl) d_{ij}| &\leq C |1-k| + C \de \sum_{j=2}^n |c_j - g_j| \, .
    \end{aligned}
\end{align}
Similarly, we define
\begin{align*}
     G = \begin{pmatrix}
        c_2 - g_2 \\ \vdots \\ c_n - g_n 
    \end{pmatrix} ,
\end{align*}
and let
\begin{align*}
    h = \frac{G \cdot \Lambda G}{2 (1-{k})^2} \, ,
\end{align*}
where $\Lambda$ is defined in \eqref{def of G, D, Lamda}. Since $\Lambda$ is positive definite, it suffices to show that $|h|$ is uniformly bounded on $\mathbb{R}$ to complete the proof. To establish this bound, we compute the derivative of $h$:
\begin{align*}
    h' = \frac{G'\cdot \Lambda G}{(1-{k})^2} + 2\frac{{k}'}{1-{k}} h \, .
\end{align*}
By \eqref{proof of S 1) k,g R}, we get
\begin{align}\label{proof of S 1) G}
    |G' - D \Lambda G| \leq C|1-k| + C \de |G| \, ,
\end{align}
and thus
\begin{align*}
    |h' - \mathbf{C}_D h| \leq C_{1/2} h^{1/2} + C_1 {\de} h + C_{3/2} {\de} |h|^{3/2} \, ,
\end{align*}
where $C_{1/2},C_1,C_{3/2}>0$ are constants depending only on $U_-$, and $\mathbf{C}_D(x)$ is a function. Since $D$ is positive definite, there exist $c_D, C_D>0$, depending only on $D,\Lambda$, such that
\begin{align*}
    c_D\leq \mathbf{C}_D \leq C_D \, .
\end{align*}
It suffices to show $h$ is bounded above by a constant on $\mathbb{R}$. We note as 
\begin{align*}
    \lim_{x\rightarrow -\infty} c_i - g_i(x) = c_i \, , \quad \lim_{x\rightarrow -\infty} 1-{k}(x) = 1 \, ,
\end{align*}
there exists $x_*\in \mathbb{R}$ and $\wt{C}>0$ such that
\begin{align*}
    h \leq \wt{C} \, , \quad on (-\infty,x_*] \, .
\end{align*}
By \eqref{proof of S 1) G}, we have
\begin{align*}
    G' \cdot \Lambda G \geq c_1 |\Lambda^{1/2} G|^2 - c_2 |\Lambda^{1/2} G||1-{k}| - c_3 {\de} |\Lambda^{1/2} G|^2 \, ,
\end{align*}
for some $c_1,c_2,c_3>0$ depending only on $U_-$.
We define $M\geq \wt{C}$ big enough such that we have \eqref{proof of S def M}, and in addition, 
\begin{align}\label{proof of S property of M R}
    c_1 - \frac{c_2}{M^{1/2}} \geq \frac{c_1}{2} \, .
\end{align}
Assume $h(x^*) = M$ for some $x^*\in \mathbb{R}$. By \eqref{proof of S property M}, we have
\begin{align*}
    h'(x^*) > \frac{c_D M}{4} \, .
\end{align*}
Hence, on $[x^*,\infty)$, we have $h\geq M$, and thus
\begin{align*}
    |\Lambda^{1/2} G| \geq M^{1/2} |1-{k}| \, .
\end{align*}
By \eqref{proof of S property of M R}, we can take $\de>0$ small enough such that on $[x^*,\infty)$,
\begin{align*}
    G' \cdot \Lambda G 
    &\geq c_1 |\Lambda^{1/2} G|^2 - \frac{c_2}{M^{1/2}} |\Lambda^{1/2} G|^2 - c_3 {\de} |\Lambda^{1/2} G|^2 \\
    &\geq \frac{c_1}{4} |\Lambda^{1/2} G|^2 \, .
\end{align*}
Hence, we get $|\Lambda^{1/2} G|\rightarrow \infty$ as $x\rightarrow \infty$. However, we know $\lim_{x\rightarrow \infty} c_i-g_i(x)=0$. We have a contradiction. Therefore, we prove that $h\leq M$ on $\mathbb{R}$.

\subsection{Proof of 2)}
By \eqref{eq for viscous shock integrated} and \eqref{Rankine-Hugoniot}, we have
\begin{align*}
    &B(S) \eta''(S) \partial_x S \cdot \mathbf{l_1} \\
    &= -\frac{(U_+-U_-)\cdot \mathbf{l_1}}{(U_+-U_-)\cdot \mathbf{l_1}} \s (S-U_-)\cdot \mathbf{l_1}
    + \frac{(U_+-U_-)\cdot \mathbf{l_1}}{(U_+-U_-)\cdot \mathbf{l_1}} \bigl(f(S)-f(U_-)\bigr) \cdot \mathbf{l_1} \\
    &= -\frac{(S-U_-)\cdot \mathbf{l_1}}{(U_+-U_-)\cdot \mathbf{l_1}} \bigl(f(U_+)-f(U_-)\bigr) \cdot \mathbf{l_1}
    + \frac{(S-U_-)\cdot \mathbf{l_1}}{(U_+-U_-)\cdot \mathbf{l_1}} \bigl(f(S)-f(U_-)\bigr) \cdot \mathbf{l_1} \\
    &\quad + \frac{(U_+-S)\cdot \mathbf{l_1}}{(U_+-U_-)\cdot \mathbf{l_1}} \bigl(f(S)-f(U_-)\bigr) \cdot \mathbf{l_1} \\
    &= \frac{(S-U_-)\cdot \mathbf{l_1}}{(U_+-U_-)\cdot \mathbf{l_1}} \bigl(f(S)-f(U_+)\bigr) \cdot \mathbf{l_1} + \frac{(U_+-S)\cdot \mathbf{l_1}}{(U_+-U_-)\cdot \mathbf{l_1}} \bigl(f(S)-f(U_-)\bigr) \cdot \mathbf{l_1} \, ,
\end{align*}
which gives
\begin{align}\label{proof of S 2) LR}
    &\frac{(U_+-U_-)\cdot \mathbf{l_1}}{(S-U_-)\cdot \mathbf{l_1} (U_+-S)\cdot \mathbf{l_1}} B(S) \eta''(S) \partial_x S \cdot \mathbf{l_1} \\
    &= \frac{\bigl(f(S)-f(U_+)\bigr) \cdot \mathbf{l_1}}{(U_+-S)\cdot \mathbf{l_1}} + \frac{\bigl(f(S)-f(U_-)\bigr) \cdot \mathbf{l_1}}{(S-U_-)\cdot \mathbf{l_1}} \\
    &=:E_R + E_L \, .
\end{align}
By Taylor expansion, \eqref{def of li}, and 1), we have
\begin{align}\label{proof of S 2) L}
    \bigl|E_L - \frac{f'(U_-)(S-U_-)\cdot \mathbf{l_1}}{\de k} - \frac{f''(U_-):(S-U_-)\otimes (S-U_-)\cdot \mathbf{l_1}}{2\de k}\bigr| &\leq C \de^2 \, ,
\end{align}
which implies
\begin{align*}
    |E_L - \lambda_1-\frac{c_f^{(1)}}{2} \de k| \leq C \de^2 \, .
\end{align*}
Similarly, we have
\begin{align*}
    \bigl|E_R - \frac{f'(U_+)(S-U_+) \cdot \mathbf{l_1}}{\de (1-k)} - \frac{f''(U_+):(S-U_+)\otimes (S-U_+) \cdot \mathbf{l_1}}{2\de (1-k)}\bigr| &\leq C \de^2 \, .
\end{align*}
Expanding further around $U_-$, we get
\begin{align*}
    \bigl|E_R - \frac{f'(U_-)(S-U_+) \cdot \mathbf{l_1}}{\de (1-k)} - \frac{f''(U_-): (U_+-U_-)\otimes (S-U_+)\cdot \mathbf{l_1}}{\de (1-k)} \\
    - \frac{f''(U_-):(S-U_+)\otimes (S-U_+) \cdot \mathbf{l_1}}{2\de (1-k)}\bigr| \leq C \de^2 \, ,
\end{align*}
which implies
\begin{align}\label{proof of S 2) R}
    |E_R + \lambda_1 + c_f^{(1)} \de - \frac{c_f^{(1)}}{2} \de (1-k)| \leq C\de^2 \, .
\end{align}
By \eqref{proof of S 2) LR}, \eqref{proof of S 2) L}, and \eqref{proof of S 2) R}, we get
\begin{align*}
    \bigl|B(S) \eta''(S) \partial_x S \cdot \mathbf{l_1} + \frac{ c_f^{(1)}}{2} \de^2 k (1-k)\bigr| \leq C \de^3 k (1-k) \, .
\end{align*}

\subsection{Proof of 3)}
By \eqref{proof of S 1) k and g} and 1), we have
\begin{equation}\label{proof of S 3) g,k}
    |g_i'-A_i^L| \leq C \de |k| \, , \quad |g_i'-A_i^R| \leq C \de |1-k| \, , \quad
    |k'| \leq C \de |k| \, , \quad |k'| \leq C \de |1-k| \, .
\end{equation}
We claim that it suffices to show that for any $2\leq i\leq n$, we have
\begin{align}\label{proof of S 3) goal}
    \bigl| \frac{\de^2 g_i'}{B(S)\eta''(S)\partial_x S \cdot \eta''(U_-)\mathbf{r_1}} \bigl| \leq C \de \, .
\end{align}
With it, \eqref{eq for viscous shock integrated each side} implies 
\begin{align*}
    \de^2 \sum_{i=2}^n|g_i'| \leq C \de^2 |k'| + C \de^3 \sum_{i=2}^n |g_i'| \, .
\end{align*}
By taking $\de>0$ small enough, we get
\begin{align*}
    \sum_{i=2}^n|g_i'| \leq C |k'| \, .
\end{align*}
Hence, it remains to prove \eqref{proof of S 3) goal}. By 2), it suffices to show that for any $2\leq i\leq n$,
\begin{align}\label{proof of S 3) goal 2}
    \bigl| \frac{g_i'}{k(1-k)} \bigl| \leq C \de \, .
\end{align}
\subsubsection{Proof of \eqref{proof of S 3) goal 2}}
Let
\begin{align*}
    G_i = \frac{(1-k^2) A_i^L + k^2 A_i^R}{k(1-k)} = \frac{1+k}{k} A_i^L + \frac{k}{1-k} A_i^R\, .
\end{align*}
By \eqref{proof of S 3) g,k}, we know 
\begin{align*}
    \bigl| \frac{g_i'}{k(1-k)} -G_i \bigl| \leq C\frac{\de |1-k^2||k| + \de|k^2||1-k|}{|k||1-k|} \leq C \de \, .
\end{align*}
Hence, to show \eqref{proof of S 3) goal 2}, it suffices to show
\begin{align}\label{proof of S 3) goal 3}
    |G_i| \leq C\de \, .
\end{align}
By \eqref{proof of S 3) A_L, A_R}, we have
\begin{align}\label{proof of S 3) def of G}
    G_i = - \frac{c_f^{(1)}}{2}d_{i1} + \sum_{j=2}^n \bigl(\lambda_j-\lambda_1\bigr) d_{ij} H_j  \, ,
\end{align}
where 
\begin{align*}
    H_j = \frac{1+k}{k} g_j - \frac{k}{1-k} (c_j-g_j) \, .
\end{align*}
We have
\begin{align*}
    H_j' = \frac{1+k}{k} g_j' + \frac{k}{1-k} g_j' + (\frac{1+k}{k})' g_j - (\frac{k}{1-k})' (c_j - g_j) \, .
\end{align*}
By \eqref{proof of S 3) g,k}, we have
\begin{align*}
    \bigl|\frac{1+k}{k} g_j' + \frac{k}{1-k} g_j' - \frac{1+k}{k} A_j^L - \frac{k}{1-k} A_j^R\bigr| &\leq C \de \, , \\
    \bigl| (\frac{1+k}{k})' g_j - (\frac{k}{1-k})' (c_j - g_j)\bigr| &\leq C \frac{|k'|}{k^2}|g_j| + C \frac{|k'|}{(1-k)^2}|c_j-g_j| \, .
\end{align*}
By \eqref{proof of S 3) g,k} and 1), we get
\begin{align*}
    \bigl|\frac{1+k}{k} g_j' + \frac{k}{1-k} g_j' - G_j\bigr| \leq C \de \, , \\
    \bigl| (\frac{1+k}{k})' g_j - (\frac{k}{1-k})' (c_j - g_j)\bigr| \leq C \de \, .
\end{align*}
Hence, we have
\begin{align*}
    |H_j' - G_j| \leq C\de \, .
\end{align*}
Together with \eqref{proof of S 3) def of G}, we get
\begin{align*}
    |G_i' - \sum_{j=2}^n \bigl(\lambda_j-\lambda_1\bigr) d_{ij} G_j | \leq C \de \, .
\end{align*}
We let 
\begin{align}
    G = \begin{pmatrix}
        G_2 \\ \vdots \\ G_n 
    \end{pmatrix}.
\end{align}
Recall $D,\Lambda$ defined in \eqref{def of G, D, Lamda}. We have 
\begin{align*}
    |G'- D \Lambda G| \leq C \de \, .
\end{align*}
We let
\begin{align*}
    h =G \cdot \Lambda G \, ,
\end{align*}
where $\Lambda$ is defined in \eqref{def of G, D, Lamda}. Since $\Lambda$ is positive definite, to complete the proof it suffices to show that $|h|$ is bounded by $\de^2$ up to a multiplicative constant, uniformly on $\mathbb{R}$. We begin by computing the derivative of $h$. We have
\begin{align*}
    |G' \cdot \Lambda G - D \Lambda G \cdot \Lambda G | \leq C_0 \de |G| \, ,
\end{align*}
which gives
\begin{align*}
    |h' - \mathbf{C}_D h | \leq C_0 \de h^{1/2} \, ,
\end{align*}
where $C_0$ is a constant depending only on $U_-$ and $\mathbf{C}_D(x)$ is a function. Since $D$ is positive definite, there exist $c_D, C_D>0$, depending only on $D,\Lambda$, such that 
\begin{align*}
    c_D \leq \mathbf{C}_D \leq C_D \, .
\end{align*}
We take $M>0$ large enough such that
\begin{align*}
    \frac{c_D M^{1/2}}{2} \geq C_0\, .
\end{align*}
Assume $h(x_*)\geq M \de^2$ at some $x_* \in \mathbb{R}$. Then
\begin{align*}
    h'(x_*) &\geq c_D h - C_0 \de h^{1/2} \\
    &\geq (c_D M^{1/2} - C_0 ) \de h^{1/2}\\
    &\geq C_0 \de h^{1/2} \, ,
\end{align*}
which implies $h\rightarrow \infty$ as $x\rightarrow \infty$. By \eqref{proof of S 3) def of G} and 1), on $\mathbb{R}$, 
\begin{align*}
    |G| \leq C \, .
\end{align*}
Hence, $h\leq C$, which gives a contradiction. We conclude that $h\leq M \de^2$ on $\mathbb{R}$.

\subsection{Proof of 4)} 
We prove by induction. We first claim that for any $j$, it suffices to show
\begin{align*}
    |B(S)\eta''(S) \partial_{x}^{j} S| \leq C \de |\partial_x S| \, ,
\end{align*}
since $B$ and $\eta''$ are positive definite. We recall \eqref{eq for viscous shock integrated}.

\noindent
\textbf{Base case $j=2$.} We have
\begin{align*}
    \bigl(B(S)\eta''(S) \partial_x S\bigr)_x = B(S)\eta''(S) \partial_{xx} S +  \bigl(B(S)\eta''(S)\bigr)_x \partial_x S \, , \\
    \bigl(-\s_1(S-U_-) + f(S) - f(U_-)\bigr)_x = \bigl(f'(S) - \s\bigr)\partial_x S \, .
\end{align*}
By \eqref{proof of S 3) g,k} and 3), we get
\begin{align*}
    |B(S)\eta''(S) \partial_{xx} S| \leq C |\partial_x S|^2 + C \de^2 |k'| + C \de^2\sum_{i=2}^n |g_i '|  \leq C \de |\partial_x S| \, .
\end{align*}

\noindent
\textbf{Induction.} Let $j_*\geq 2$. Assume that the inequality holds for all $2\leq j\leq j_*$. We then apply $\partial_x^{j_*}$ to both sides of the first equation in \eqref{eq for viscous shock integrated}. We have
\begin{align*}
    |B(S)\eta''(S) \partial_{x}^{j_*+1} S| \leq C |\partial_x S|^2 + C |\partial_{x}^{j_*} S| \leq C \de |\partial_x S| \, .
\end{align*}

\medskip
\medskip

\noindent
\textbf{Declaration on the use of AI.} During the preparation of this manuscript, ChatGPT (OpenAI) was used to assist with language editing, checking mathematical arguments for potential errors, and checking references. The authors reviewed and verified all resulting revisions and take full responsibility for the content of the manuscript.

\bibliographystyle{abbrv}
\bibliography{refs}

\end{document}